\documentclass[a4paper, 11pt]{article}
\usepackage[utf8]{inputenc}
\usepackage[T1]{fontenc}
\usepackage{tgpagella}
\usepackage{amsmath, amssymb, amsthm, bm, xcolor}
\usepackage{algorithm}
\usepackage{algpseudocode}
\usepackage{enumitem}

\usepackage[textwidth = 15.2cm]{geometry}
\definecolor{main_color}{HTML}{2789d6}
\usepackage[pdfa, hidelinks, colorlinks = true, allcolors = main_color]{hyperref}

\newtheorem{theorem}{Theorem}
\newtheorem{lemma}{Lemma}
\newtheorem{proposition}{Proposition}
\newtheorem{corollary}{Corollary}

\DeclareMathOperator{\conv}{conv}

\newcommand{\R}{\mathbb{R}}
\newcommand{\Z}{\mathbb{Z}}
\title{Sharp Proximity Bounds for Integer Programs with Few Rows}
\author{Stefan Kuhlmann\thanks{ETH Zürich, Department of Mathematics, Institute for Operations Research, Zürich, Switzerland, \texttt{stefan.kuhlmann@math.ethz.ch}} \and Luze Xu\thanks{Hong Kong University of Science and Technology, Department of Mathematics, Hong Kong, China, \texttt{xuluze@ust.hk}}}
\date{ }
\begin{document}
\maketitle

\begin{abstract}
\noindent
We study the proximity between an optimal vertex of a standard form linear program with constraint matrix $A\in\Z^{m\times n}$ and an optimal integer solution. When $m = 2$, we prove the sharp infinity-norm bound $\Delta-1$, where $\Delta$ is the largest absolute maximal minor of $A$. When $m = 1$, we express the difference between an optimal LP vertex and an optimal integer solution as a nonnegative combination of certain generators in the kernel of $A$ with a bounded coefficient sum. This implies sharp proximity bounds for all norms when $m=1$. Combining this representation with its known counterpart for Hilbert bases due to Sissokho, we develop techniques that transfer results from $m=1$ to $m=2$. Along the way, for $m=2$, we prove that such a representation exists for an extended Gomory group relaxation and obtain the sharp infinity-norm bound $\Delta$ for Graver basis elements. We also present a new short proof of Sissokho's result and lower bounds on the proximity for all $m$ and $\Delta$.
\end{abstract}

\section{Introduction}\label{sec:introduction}
Proximity bounds estimate the distance from an optimal vertex of a linear relaxation of an integer program to a nearest optimal integer solution. They connect continuous and discrete optimization and play a role in sensitivity analysis~\cite{cook1986sensitivity} and algorithm design~\cite{aprileIntegerProgramsNearly2024,eisenbrand_proximity_2020,hochbaumConvexSeparableOptimization1990}. 
The classical proximity result of Cook, Gerards, Schrijver, and Tardos~\cite{cook1986sensitivity} bounds the infinity-norm proximity in terms of the number of variables and subdeterminants of the constraint matrix. Their bound is known to be tight for rational right-hand sides \cite{schrijver1998theory}. For integer-valued right-hand sides, recent work has improved this bound and extended the improvement to arbitrary norms \cite{celaya2022proximity}. Determining sharp proximity bounds for general integer programs remains an open problem. 
In this paper, we study sharp proximity bounds for standard form integer programs with one or two equations.

For $A\in\Z^{m\times n}$ of full row rank, $b\in\Z^m$, and $c\in\R^n$, consider
\begin{align}
\mathrm{LP}(A,b,c)&=\max\{c^\top x:Ax=b,\ x\in\R_{\ge0}^n\},\label{eq:LP}\tag{LP}\\
\mathrm{IP}(A,b,c)&=\max\{c^\top z:Az=b,\ z\in\Z_{\ge0}^n\}.\label{eq:IP}\tag{IP}
\end{align}
Assume that \eqref{eq:IP} attains a finite optimum, which implies that \eqref{eq:LP} has optimal vertex solutions.  
For a column submatrix $A_I$ with column indices in $I$, put 
$\Delta(A)=\max_{|I|=m}\left|\det A_I\right|$ and $\gcd(A)=\gcd\{\left|\det A_I\right|:|I|=m\}$.
We assume $\gcd(A)=1$ throughout the paper. This can be done without loss of generality because left multiplication of $(A,b)$ by the inverse of a basis matrix of the lattice $A\Z^n$ preserves the feasible points and gives integral data with $\gcd(A)=1$. 
For a fixed instance $(A,b,c)$, let $V^*(A,b,c)$ be the set of optimal LP vertices and $Z^*(A,b,c)$ the set of optimal integer solutions. Let $\Vert\cdot\Vert_p$ denote the $\ell_p$-norm for $1\leq p\leq \infty$ and define 
\[
 \pi^p(A,b,c)=\max_{x^*\in V^*(A,b,c)}\, \min_{z^*\in Z^*(A,b,c)}\|z^*-x^*\|_p. 
\]
For positive integers $m$ and $\Delta$, define 
\[
 \pi^p(m,\Delta)=\sup_{A,b,c}\pi^p(A,b,c), 
\]
where the supremum ranges over all $n\ge m$, all matrices $A\in\Z^{m\times n}$ of full row rank with $\gcd(A)=1$ and $\Delta(A)=\Delta$, and all $b\in\Z^m$, $c\in\R^n$ for which \eqref{eq:IP} is feasible and attains a finite optimum.

It is indeed possible to bound $\pi^p(m,\Delta)$ solely by $m$ and $\Delta$. A change of norm in a proximity result due to Eisenbrand and Weismantel \cite{eisenbrand_proximity_2020} implies $\pi^p(m,\Delta)\leq \pi^1(m,\Delta)\leq m(2m+1)^m\Delta$. For fixed $\Delta$, one can improve the dependence on $m$ by using sparsity. Following the approach of~\cite{lee2020improving}, we combine Cook et al.'s proximity bound~\cite{cook1986sensitivity} with the sharp sparsity bound of $m+\Delta-1$~\cite[Theorem~4]{kuhlmannThresholdPhenomenonShortest2026} to obtain 
\begin{equation*}
    \pi^p(m,\Delta)\leq (m+1)^{1/p}(m+\Delta-1)\Delta.
\end{equation*}
One can improve this bound by a factor of approximately $4/9$ using \cite{celaya2022proximity} and, for varying $\Delta$, adapt it using other available sparsity bounds \cite{alievSupportIntegerOptimal2018,gribanov$$Delta$$modularInteger2024}. None of these bounds are believed to be sharp. For instance, for the $\ell_{\infty}$-norm, generalizing the best-known lower bound \cite{aliev2020distances} combined with the upper bound above imply
\begin{equation*}
\Delta - 1\leq \pi^\infty(m,\Delta)\leq ((4(m+\Delta-1)+2)/9) \Delta.
\end{equation*}
It is open whether $\pi^\infty(m,\Delta)$ depends on $m$. The few known tight bounds suggest that it does not: one has $\pi^\infty(m,1) = 0$, $\pi^\infty(m,2) = 1$ \cite[Theorem~6]{celaya2022proximity}, and computer enumerations have confirmed that the lower bound is tight for bounded feasible regions when $(m,\Delta)=(2,3), (2,4), (2,5), (3,3)$ \cite{koppe2026determining}. 
If one does not require optimality of the integer solution, i.e., restricts to $c=0$, Aliev, Henk, and Oertel have shown that the bound $\Delta - 1$ is tight for $m = 1$ \cite{aliev2020distances}. Our first result determines a sharp proximity bound for optimal solutions when $m = 2$. 

\begin{theorem}\label{thm:infinity}
For every positive integer $\Delta$, $\pi^\infty(2,\Delta)=\Delta-1$.
\end{theorem}
Thus, for every optimal LP vertex $x^*$ of a two-row instance, there exists an optimal integer solution $z^*$ with $\|z^*-x^*\|_\infty\le\Delta-1$, and equality is attained for every positive integer $\Delta$. Our approach for proving the sharp bound in Theorem~\ref{thm:infinity} is to prove a structural result when $m = 1$ that can be lifted to $m = 2$.

Let $e_i$ denote the $i$-th standard unit vector. 
For $A=a^\top$, since $b = 0$ implies $x^*\in \Z^n$ and thus $z^*-x^* = 0$, we may assume without loss of generality that $a_1,b>0$ and $x^*=(b/a_1)e_1$. 
Let $\Delta=\max_{i\in[n]}|a_i|$ and partition the nonbasic set $N=\{2,\ldots,n\}$ into $I_+=\{i\in N:a_i>0\}$, $I_-=\{j\in N:a_j<0\}$, and $I_0=\{i\in N:a_i=0\}$. 
Define the normalized generators 
\begin{align*}
 r_{1i}&=-a_i e_1+a_1e_i &&(i\in I_+\cup I_-),\\ r_{ij}&=-a_j e_i+a_ie_j &&(i\in I_+,\ j\in I_-), 
\end{align*}
and let $\mathcal R(a)$ be their collection. These are integer vectors in $\ker(a^\top)$ with exactly two nonzero coordinates and nonnegative nonbasic coordinates. Note that we do not divide these generators by the greatest common divisors of their coordinates. A vector $u$ is conformal to $v$ if $u_iv_i\ge0$ for all $i$ with $v_i\neq 0$ and $u_i=0$ whenever $v_i=0$.
\begin{theorem}\label{thm:main}
For every optimal LP vertex $x^*$ with the basis $\{1\}$, there exists an optimal integer solution $z^*$ such that 
\[
 z^*-x^*\in\left(1-\frac1{\Delta}\right) \conv(\{0\}\cup\mathcal R(a)). 
\]
Every generator used in the representation can be chosen conformal to $z^*-x^*$.
\end{theorem}
It is possible to prove a refined version of Theorem~\ref{thm:main} with $1-1/\Delta$ replaced by $1-1/a_1$. Since we do not need this refined statement, we omit it. 

By applying the triangle inequality and showing that one can exclude $r_{ij}\in \mathcal{R}(a)$ with $\gcd(r_{ij}) = \Delta$ in Theorem~\ref{thm:main}, we get sharp proximity bounds when $m=1$ for any norm, in particular, $\ell_p$-norms.
\begin{corollary}\label{cor:allp}
For every positive integer $\Delta$, 
\[
 \pi^\infty(1,\Delta)=\Delta-1,\qquad \pi^p(1,\Delta)=\left(1-\frac1\Delta\right) \bigl(\Delta^p+(\Delta-1)^p\bigr)^{1/p}\quad(1\le p<\infty). 
\]
In particular, 
\[
 \pi^1(1,\Delta)=2\Delta-3+1/\Delta,\qquad \pi^2(1,\Delta)=\left(1-\frac1\Delta\right) \sqrt{\Delta^2+(\Delta-1)^2}. 
\]
\end{corollary}
Section~\ref{sec:lower-bound} provides the matching lower bound examples for Theorem~\ref{thm:main} and Corollary~\ref{cor:allp}. Their feasible regions are rays of the form $z^* = x^* + \lambda r_{ij}$ for $\lambda\geq 0$ and some $r_{ij}$ maximizing the norm in question over all normalized generators. This low-dimensionality is related to research that investigates the trade-off between proximity and sparsity of an optimal integer solution \cite{alievSparsityProximityTransference2025,alievDistanceSparsityTransferenceVertices2021}. In \cite[Theorem~1]{alievSparsityProximityTransference2025}, Aliev, Celaya, and Henk show that an upper bound on the integrality gap proportional to $\sqrt{\det (AA^\top)}$ gets exponentially smaller for larger support of $z^*$. They translate this result to one in terms of $m$ and $\Delta$ with the same asymptotic dependence on the support of $z^*$ \cite[Corollary~2]{alievSparsityProximityTransference2025}. These results suggest further refinements of Theorem~\ref{thm:main}.

To the best of our knowledge, the estimates in Corollary~\ref{cor:allp} are the first sharp results for optimal integer solutions when $m=1$ and general $a^\top$. 
There exists an asymptotically tight proximity bound for the Euclidean norm: generalizing a feasibility result \cite[Theorem~3]{alievSparsityProximityTransference2025} to optimal integer solutions, Aliev, Averkov, Jones, and Oertel \cite[Theorem~1]{aliev2026optimal} prove an asymptotically tight upper bound of $\sqrt{\det(AA^\top)} - 1$ for all $m$. When $m=1$, this equals $\Vert a\Vert_2-1$. Our sharp $\ell_2$-norm estimate in Theorem~\ref{thm:main}, after applying the triangle inequality, implies a stronger bound depending only on the two largest coefficients of $a$ in absolute value. This can be upper bounded by $(1-1/\Delta)\|a\|_2$. Corollary~\ref{cor:allp} gives the sharp bound that only depends on $\Delta$.

The convex hull result in Theorem~\ref{thm:main} is motivated by a result due to Sissokho for Hilbert basis elements \cite{sissokhoGeometryMinimalSolutions2021}. For an integer matrix $A$, the \emph{Hilbert basis} $\mathcal H(A)$ consists of the nonzero vectors in $\ker(A)\cap\Z_{\ge0}^n$ that cannot be written as a sum of two nonzero vectors in this set. When $m=1$, Sissokho shows that every Hilbert basis element lies in the convex hull of zero and the nonnegative normalized generators.
\begin{theorem}[{\cite[Theorem~1.2]{sissokhoGeometryMinimalSolutions2021}}]\label{thm:hilbert}
Let $a\in \Z^n$ have no zero coordinate, $J_{+} = \lbrace i \in \lbrack n\rbrack : a_i > 0\rbrace$ and $J_{-} = \lbrace j \in \lbrack n\rbrack : a_j < 0\rbrace$. For $h\in\mathcal H(a^\top)$, 
\[
 h\in\conv\left(\{0\}\cup\{r_{ij}:i\in J_+,\ j\in J_-\}\right). 
\]
\end{theorem}
Theorem~\ref{thm:hilbert} affirmatively answered a conjecture by Henk and Weismantel \cite[Conjecture~1]{HenkWeismantel2000} and, independently, by Hoşten and Sturmfels \cite[p. 54]{HenkWeismantel2000}. In Section~\ref{sec:one-proof}, we present a new short proof of Theorem~\ref{thm:hilbert}. 

The advantage of having statements in terms of the normalized generators for $m$ rows is that one can lift them to statements for $m + 1$ rows. We show in Lemma~\ref{lem:two-row-bases} that Theorem~\ref{thm:hilbert} already implies refined $\ell_\infty$-norm bounds of the form $\Vert h\Vert_\infty\leq \Delta$ for Hilbert basis elements. Then, using Lemma~\ref{lem:two-row-bases} and a rounding technique that relates $z^*-x^*$ to a Hilbert basis element, already gives an almost tight upper bound of $\Delta - 1/2$ on the $\ell_\infty$-proximity; see the nonbasic variable part of the proof of Theorem~\ref{thm:infinity}. 
For our proximity variant, Theorem~\ref{thm:main} lifts to a convex hull statement for extended Gomory relaxations; see Theorem~\ref{thm:extended-gomory}. 
This result already implies a variant of Theorem~\ref{thm:main} when $m = 2$ for all but finitely many right-hand sides for each fixed matrix $A$. Combining both approaches proves Theorem~\ref{thm:infinity}.

Since the lifting strategy works for general $m$, the natural next step is to show variants of Theorem~\ref{thm:main} and Theorem~\ref{thm:hilbert} when $m=2$. This is likely more challenging: Theorem~\ref{thm:main} and Theorem~\ref{thm:hilbert} hold for conformal normalized generators. It was conjectured that the conformal representation in Theorem~\ref{thm:hilbert} extends to arbitrary $m$ \cite[cf. Strong Hilbert basis conjecture]{HenkKuhlmannWeismantel2022}. We refute this conjecture by showing that the normalized generators cannot always be chosen to be conformal when $m=2$; see Appendix~\ref{app:conformal}. The weaker conjecture \cite[p. 1929]{HenkKuhlmannWeismantel2022} remains open and is equivalent to Theorem~\ref{thm:hilbert} for general $m$ when allowing nonconformal normalized generators. Finally, let us remark that using nonconformal normalized generators does not cause an issue with the lifting strategy proposed in this paper.

In Section~\ref{sec:link-hilbert-bases}, we introduce our notation and establish the link between proximity and Hilbert bases. Section~\ref{sec:one-proof} proves our results when $m=1$. Section~\ref{sec:two-ingredients} demonstrates how to lift Theorem~\ref{thm:hilbert} to results on Hilbert and Graver basis elements that hold for two rows. Section~\ref{sec:two-proof} explains the lifting of Theorem~\ref{thm:main} and completes the proof of Theorem~\ref{thm:infinity}. Section~\ref{sec:lower-bound} gives the matching lower-bound examples.

\section{Notation and Linking Proximity with Hilbert Bases}\label{sec:link-hilbert-bases}
We write $[n]=\{1,\ldots,n\}$. For any $i,j\in[n]$, $a_i$ denotes column $i$ of $A$, and $A_{-j}$ is $A$ with column $j$ deleted. For $I\subseteq[n]$, $x_I$ denotes the subvector of $x$ indexed by $I$, and, for a vector $v$, $|v|$ denotes its componentwise absolute value. When the matrix $A$ is fixed, we abbreviate $\Delta:=\Delta(A)$. Recall that a vector $u$ is conformal to $v$ if $u_iv_i\ge0$ for all $i$ with $v_i\neq 0$ and $u_i=0$ whenever $v_i=0$. The \emph{Graver basis} $\mathcal{G}(A)$ consists of the nonzero integer kernel vectors $g$ that cannot be written as $g=u+v$ with nonzero integer kernel vectors $u,v$ both conformal to $g$. Every nonnegative integer kernel vector decomposes into Hilbert basis elements, and every integer kernel vector decomposes into Graver basis elements conformal to it. For a diagonal sign matrix $S$ and $h\in\Z_{\ge0}^n$, $h\in\mathcal H(AS)$ if and only if $Sh\in\mathcal G(A)$. So changing column signs preserves absolute minors and sends a Graver element $g$ to the Hilbert basis element $h=|g|$ of the sign-adjusted matrix. It will therefore be enough to prove both bounds for $h\in\mathcal H(A)$ and the statements in this paper remain valid for Graver basis elements.

For a full-row-rank matrix $A$, a basis $B$, and $N=[n]\setminus B$, an integer feasible point $z$ is minimal in its nonbasic coordinates if 
\begin{equation}
\tag{$\ast$}\label{eq:M} \nexists\,w\in\Z_{\ge0}^n: Aw=b,\quad w_N\le z_N,\quad w_N\ne z_N. 
\end{equation}
For an optimal LP vertex $x^*$ with a corresponding optimal basis $B$, we have 
\[
 x_B^*=A_B^{-1}b,\qquad x_N^*=0,\qquad \bar c_N^\top=c_N^\top-c_B^\top A_B^{-1}A_N\le0. 
\]
Every feasible $x$ satisfies $c^\top x=c^\top x^*+\bar c_N^\top x_N$. Choose $z^*\in Z^*(A,b,c)$ minimizing $\sum_{i\in N}z_i^*$. A point violating \eqref{eq:M} would also be optimal and would have a smaller nonbasic sum.
Thus this selected $z^*$ satisfies \eqref{eq:M}. It will provide a proof of our proximity results.

We next establish the connection between proximity and Hilbert bases.
\begin{lemma}\label{lem:augmented}
Let $A\in\Z^{m\times n}$ have full row rank, let $B$ be a basis, and let $z\in\Z_{\ge0}^n$ satisfy $Az=b$. 
If Condition~\eqref{eq:M} holds, then, for every $\bar x\in\Z_{\ge0}^n$ supported on $B$, 
\[
 (|z-\bar x|,1)\in\mathcal H([AS,-(b-A\bar x)]), 
\]
where $S$ is a diagonal matrix with $S_{ii}=\operatorname{sign}(z_i-\bar x_i)$ when $z_i-\bar x_i\neq 0$, and $S_{ii}\in \lbrace -1,1\rbrace$ otherwise.
\end{lemma}
\begin{proof}
It is easy to verify that $(|z-\bar x|,1)\in \ker([AS,-(b-A\bar x)])\cap\Z_{\ge 0}^{n+1}$. Suppose that $(|z-\bar x|,1)\notin \mathcal H([AS,-(b-A\bar x)])$. Then there exists a nonzero kernel subvector with final coordinate zero, denoted by $(u,0)$ with $0\le u\le|z-\bar x|$ and $ASu=0$. Then $w = z - Su$ satisfies $Aw = Az - ASu =b$.
Also, $w\ge 0$ because it lies componentwise between $z$ and $\bar x$. Moreover, $w_N = z_N - u_N\le z_N$ and $w_N\ne z_N$, because otherwise $u_N=0$ and the nonsingularity of $A_B$ would imply $u=0$. Thus $w$ violates \eqref{eq:M}.
\end{proof}

\section{Proofs of the One-Row Results}\label{sec:one-proof}
We first present an alternative and shorter proof of Theorem~\ref{thm:hilbert}. We then apply the same tools that provide a proof of Theorem~\ref{thm:hilbert} to prove Theorem~\ref{thm:main}. Suppose that $a$ has no zero coordinates. 
For some element $v\in \ker(a^\top)$, there are usually many choices of normalized generators to express $v$ as a nonnegative combination. Even if $v$ is contained in the convex hull of some of those generators, it could still be that several nonnegative combinations are not convex combinations. We claim that the following choices of coefficients yield proofs of Theorem~\ref{thm:hilbert} and Theorem~\ref{thm:main}. Let $J_+=\{i:a_i>0\}$ and $J_-=\{j:a_j<0\}$. If $v\geq 0$ and $v\neq 0$, then one can express $v$ as
\begin{equation}\label{eq:convex-positive}
    v = \sum_{i\in J_+, j\in J_-} \frac{v_iv_j}{\sum_{j\in J_-}v_j(-a_j)} r_{ij}.
\end{equation}
For $v_1 < 0$ and $v_i\geq 0$ for $i\geq 2$, set $I_+ = J_+ \backslash \lbrace 1 \rbrace$ and $I_- = J_-$. Then one gets 
\begin{equation}\label{eq:convex-negative}
    v = \sum_{i\in I_+, j\in I_-} \frac{v_iv_j}{\sum_{i\in I_+}v_ia_i} r_{ij} + \sum_{i\in I_+} \frac{v_i(-v_1)}{\sum_{i\in I_+}v_ia_i}r_{1i}.
\end{equation}
Each coefficient is nonnegative, and one can verify both equations by plugging in the definitions of the generators. The challenge in proving Theorem~\ref{thm:hilbert} and Theorem~\ref{thm:main} is to show that the coefficient sum is indeed small. To do so, we use another result due to Sissokho.
\begin{lemma}[Theorem 1.3 in \cite{sissokho2014note}]\label{lem:product}
Let $h\in\mathcal H(a^\top)$ with $J_+$ and $J_-$ from above. Then 
\[
    \left(\sum_{i\in J_+}h_i\right)\left(\sum_{j\in J_-}h_j\right)\le \sum_{i\in J_+}a_i h_i = -\sum_{j\in J_-}a_j h_j
\]
and 
\[
 \sum_{i\in J_+}h_i\le\max_{j\in J_-,\,h_j>0}(-a_j),\qquad \sum_{j\in J_-}h_j\le\max_{i\in J_+,\,h_i>0}a_i. 
\]
\end{lemma}

The proof of Theorem~\ref{thm:hilbert} is now merely an application of Lemma~\ref{lem:product}.
\begin{proof}[Proof of Theorem~\ref{thm:hilbert}]
Use the nonnegative combination in \eqref{eq:convex-positive} to express $h$. The coefficient sum equals
\begin{equation*}
     \sum_{i\in J_+, j\in J_-} \frac{h_ih_j}{\sum_{j\in J_-}h_j(-a_j)} = \frac{\left(\sum_{i\in J_+} h_i\right)\left(\sum_{j\in J_-} h_j\right)}{\sum_{j\in J_-}h_j(-a_j)}\leq 1,
\end{equation*}
where the inequality follows from Lemma~\ref{lem:product}.
\end{proof}

We use the following well-known estimate on the Davenport constant of finite abelian groups. Its proof follows from the pigeonhole principle.
\begin{lemma}\label{lem:zero-sum}
Let $G$ be a finite abelian group of order $k$. If a sequence $g_1,\ldots,g_s\in G$ has no nonempty subsequence with sum zero, then $s\le k-1$.
\end{lemma}
\begin{proof} 
If $s\ge k$, two of the $k+1$ partial sums $0,g_1,g_1+g_2,\ldots,g_1+\cdots+g_k$ are equal. Their difference is a nonempty consecutive subsequence with sum zero, a contradiction.
\end{proof}
\begin{proof}[Proof of Theorem~\ref{thm:main}]
Choose $z^*\in Z^*(A,b,c)$ minimizing $\sum_{i\in N}z_i^*$. Then $z^*$ satisfies \eqref{eq:M}. Notice that \eqref{eq:M} forces $z^*_i = 0$ if $a_i = 0$. So we assume that $a$ has no zero coordinates. 
If $x^*$ is integral, then $x_N^*=0$, so \eqref{eq:M} forces $z^*_N=0$ and hence $z^*=x^*$. Assume now that $x^*$ is fractional, and write $d=z^*-x^*\in \ker(a^\top)$. Then $d_1\neq 0$ is fractional and, depending on whether $d_1 > 0$ or $d_1 < 0$, we will show that either \eqref{eq:convex-positive} or \eqref{eq:convex-negative} has the desired coefficient sum. Before that we fix some notation. 
Let $N=\{2,\ldots,n\}$ and $I_+=\{i\in N:a_i>0\}$, $I_-=\{j\in N:a_j<0\}$. 
Write $x^*_1 = \lfloor x^*_1\rfloor + r/a_1$ for $1\le r<a_1$ and, for lighter calculations, define 
\[
 H=\sum_{i\in I_+}z^*_i,\quad K=\sum_{j\in I_-}z^*_j,\quad P=\sum_{i\in I_+}a_i z^*_i,\quad Q=\sum_{j\in I_-}(-a_j)z^*_j. 
\]
We split the proof according to whether $d_1<0$ or $d_1 > 0$. 

\emph{Case 1: $d_1 < 0$.} Let $d$ be expressed as in \eqref{eq:convex-negative} and $\ell = \lfloor x^*_1\rfloor - z^*_1\geq 0$. We want to estimate the coefficient sum
\begin{equation}
    \label{eq:prox-coeffsum-negative}
    \sum_{i\in I_+, j\in I_-} \frac{z^*_iz^*_j}{\sum_{i\in I_+}z^*_ia_i} + \sum_{i\in I_+} \frac{z^*_i(x^*_1 - z^*_1)}{\sum_{i\in I_+}z^*_ia_i} = \frac{H(K + \ell + r/a_1)}{P}.
\end{equation}
For that purpose, we derive additional inequalities. For $s = \lfloor x^*_1\rfloor$ and $s = \lceil x^*_1 \rceil$, apply Lemma~\ref{lem:augmented} with $\bar{x}=se_1$ to obtain $ (|z^*_1-s|,z^*_N,1)\in \mathcal H\bigl((-a_1,a_N^\top,-(b-a_1s))\bigr)$. In the case $s = \lfloor x^*_1 \rfloor$, the positive and negative multiplicities and their common weighted sum are $(H,K+\ell+1,P)$. In the other case, we get $(H+1,K+\ell+1,P+a_1-r)$. Thus, by Lemma~\ref{lem:product}, 
\[
 H(K+\ell+1)\le P,\qquad (H+1)(K+\ell+1)\le P+a_1-r. 
\]
Suppose towards a contradiction that $H\ge a_1$. Apply Lemma~\ref{lem:zero-sum} to the $H$ occurrences of the coefficients $a_i$, $i\in I_+$, in $\Z/a_1\Z$. There are integers $0\le w_i\le z^*_i$, not all zero, such that $\sum_{i\in I_+}a_iw_i\equiv0\pmod{a_1}$. Set $z'_1=z^*_1+\sum_{i\in I_+}a_iw_i/a_1$, $z'_i=z^*_i-w_i$ for $i\in I_+$, and $z'_i=z^*_i$ for $i\in N\setminus I_+$. Then $z'$ is an integer feasible point with a strictly smaller nonbasic part, contradicting \eqref{eq:M}. This establishes $H\le a_1-1$. Hence, we can multiply the first of the two inequalities above by $a_1-1-H\ge0$ and the second by $H\geq 0$. Adding them gives $a_1H(K+\ell+1)\le(a_1-1)P+H(a_1-r)$. Rearranging terms gives 
\begin{equation*}
    \frac{H(K+\ell+r/a_1)}P\le1-\frac{1}{a_1}\leq 1 - \frac{1}{\Delta}.
\end{equation*}
Using this inequality in \eqref{eq:prox-coeffsum-negative} proves the claim when $d_1 < 0$.

\emph{Case 2: $d_1>0$.} Let $d$ be expressed as in \eqref{eq:convex-positive} and $u = z^*_1-\lfloor x^*_1\rfloor \geq 1$. We want to estimate the coefficient sum
\begin{equation}
    \label{eq:prox-coeffsum-positive}
    \sum_{i\in I_+, j\in I_-} \frac{z^*_iz^*_j}{\sum_{j\in I_-}z^*_j(-a_j)} + \sum_{j\in I_-} \frac{z^*_j(z^*_1 - x^*_1)}{\sum_{j\in I_-}z^*_j(-a_j)} = \frac{K(H + u - r/a_1)}{Q}.
\end{equation}
As in the previous case, we derive additional inequalities. For $s = \lfloor x^*_1\rfloor$ and $s = \lceil x^*_1 \rceil$, apply Lemma~\ref{lem:augmented} with $\bar{x}=se_1$ to obtain $ (|z^*_1-s|,z^*_N,1)\in \mathcal H\bigl((a_1,a_N^\top,-(b-a_1s))\bigr)$. In the case $s = \lfloor x^*_1 \rfloor$, the positive and negative multiplicities and their common weighted sum are $(H+u,K+1,Q+r)$. In the other case, we get $(H+u,K,Q)$. Lemma~\ref{lem:product} gives 
\[
 (H+u)K\le Q,\qquad (H+u)(K+1)\le Q+r. 
\]
Next, we need to show that $K\le \Delta - 1$. To derive this inequality, apply Lemma~\ref{lem:augmented} with $\bar{x} = 0$ to obtain $(z^*,1)\in \mathcal{H}([a^\top,-b])$. Applying Lemma~\ref{lem:product} to $(z^*,1)$ gives $K + 1\le \max_{i \in I_+\cup \lbrace 1\rbrace} a_i \leq \Delta$. So $K\leq \Delta - 1$. Here $Q>P\ge0$, so $K\ge1$. Hence, we can multiply the first inequality above by $\Delta - 1 - K\geq 0$ and the second by $K$. Adding them gives $\Delta(H+u)K\le (\Delta-1)Q + rK$. Since $\Delta\ge a_1$, this yields
\[
 \Delta K(H+u-r/a_1)\le (\Delta-1)Q+rK(1-\Delta/a_1)\le (\Delta-1)Q,
\]
and therefore
\begin{equation*}
    \frac{K(H+u-r/a_1)}Q\le1-\frac{1}{\Delta}.
\end{equation*}
Using this in \eqref{eq:prox-coeffsum-positive} settles the case $d_1 > 0$.
\end{proof}
\begin{proof}[Proof of Corollary~\ref{cor:allp}]
We claim that the representation in Theorem~\ref{thm:main} can avoid every generator whose two nonzero coordinates both have absolute value $\Delta$. If so, every used generator has one entry of magnitude at most $\Delta$ and the other of magnitude at most $\Delta-1$, and the claim follows from the triangle inequality. So it remains to rule out that a used generator has two coordinates with absolute value $\Delta$. There are three possibilities.
If $a_i=\Delta$ and $a_j=-\Delta$ for some $i\in I_+$ and $j\in I_-$, the pair $r_{ij}$ cannot be used, because deleting $e_i+e_j$ would violate \eqref{eq:M}. If $a_1=a_i=\Delta$ for some $i\in I_+$, the generator $r_{1i}$ cannot be used, because deleting $e_i$ and adding $e_1$ would violate \eqref{eq:M}.
Finally, if $a_1=-a_j=\Delta$ for some $j\in I_-$ and $z_j^*>0$, minimality forces $z_1^*=0$, because otherwise deleting $e_1+e_j$ would violate \eqref{eq:M}. Thus $d_1<0$, and the representation \eqref{eq:convex-negative} in Case 1 uses no such $r_{1j}$. This finishes the proof of Corollary~\ref{cor:allp}.
\end{proof}

\section{Proofs of the Two-Row Hilbert and Graver Bounds}\label{sec:two-ingredients}
We first put a rank-two matrix $A\in\Z^{2\times n}$ into a standard form that will also be used in Section~\ref{sec:two-proof}. Fix a nonzero column $a_j$, and let $k>0$ be the greatest common divisor of its coordinates. There exists a unimodular matrix $U\in\Z^{2\times2}$ such that $Ua_j=(k,0)^\top$.
Left multiplication by $U$ preserves the kernel, the Hilbert and Graver bases, and all absolute $2\times2$ minors. Replacing $(A,b)$ by $(UA,Ub)$ also preserves all feasible points. We keep the notation $A,b$ after this transformation and, displaying column $j$ first, write 
\begin{equation}
\label{eq:two-row-standard} A=\begin{pmatrix}k&p^\top\\0&q^\top\end{pmatrix},\qquad b=\binom{\gamma}{\beta},\qquad a_i=\begin{pmatrix} p_i \\ q_i\end{pmatrix}\quad(i\ne j). 
\end{equation}
The vectors $p,q$ are indexed by $[n]\setminus\{j\}$, and $q\ne0$ because $A$ has rank two. In this form, 
\begin{equation}
    \label{eq:two-row-determinants}
 \det(a_j,a_i)=kq_i,\qquad \det(a_i,a_\ell)=p_iq_\ell-p_\ell q_i\quad(i,\ell\ne j).
\end{equation}
Notice that, if $q_i = 0$ and $u\in\mathcal{H}(q^\top)$ with $u_i > 0$, then $u = e_i$. 
If $B=\{j,\ell\}$ is a basis, the sign of the second row can be chosen so that $q_\ell=t>0$.
Then $kt=\left|\det A_B\right|$ and, if $x^*$ is represented by this basis, $\beta=t x_\ell^*$. 
We use the normalization from above to prove the following.
\begin{lemma}\label{lem:two-row-bases}
Let $A\in\Z^{2\times n}$ have rank two.
\begin{enumerate}[label=(\roman*),ref=(\roman*)]
\item\label{part:strict-column} Suppose that the absolute values of all $2\times 2$ minors are at most a positive integer $\Delta$ and that $|\det(a_e,a_i)|<\Delta$ for every $i\in[n]$. Then every $h\in\mathcal H(A)$ with $h_e>0$ satisfies $h_j<\Delta$ for $j\ne e$.
\item\label{part:graver-coordinate} For every $g\in\mathcal G(A)$ and every $j\in[n]$, $|g_j|\le\Delta(A_{-j})$.
\end{enumerate}
\end{lemma} 
\begin{proof} By the discussion in Section~\ref{sec:link-hilbert-bases}, it suffices to prove the statement for Hilbert basis elements.
Let us first describe the lifting that provides a proof of \ref{part:strict-column} and \ref{part:graver-coordinate}.
Fix $j$ with $h_j>0$. If $a_j=0$, then $h=e_j$, so $h_j=1\le\Delta(A_{-j})$ by rank two. For part~\ref{part:strict-column}, we have $j=e$ and $h_k=0<\Delta$ for all $k\ne e$. So we may assume that $a_j\ne0$. Use \eqref{eq:two-row-standard}. Because $q^\top h_{-j}=0$, decompose 
\[
 h_{-j}=u^1+\cdots+u^s,\qquad u^r\in\mathcal H(q^\top). 
\]
Then the first row implies $kh_j+\sum_{r=1}^s p^\top u^r=0$. We claim that this equation gives rise to another Hilbert basis element, namely $(h_j,1,\ldots,1) \in \mathcal H((k,p^\top u^1,\ldots,p^\top u^s))$. 
Suppose that is not the case. Then there exists a nonzero proper nonnegative integer kernel subvector $(t,\varepsilon_1,\ldots,\varepsilon_s)$ of $(h_j,1,\ldots,1)$, with $0\le t\le h_j$, $\varepsilon_r\in\{0,1\}$, and $kt+\sum_{r=1}^s\varepsilon_r p^\top u^r=0$. Define $y\in \Z^n_{\geq 0}$ by $y_j = t$ and $y_{-j}= \sum_{r=1}^{s} \varepsilon_ru^r$. We have $0\le y\le h$, $y\in \ker(A)$, $y\neq 0$, and $y \neq h$ by construction. This contradicts $h\in\mathcal H(A)$. So we must have $(h_j,1,\ldots,1) \in \mathcal H((k,p^\top u^1,\ldots,p^\top u^s))$. In particular, $p^\top u^r\ne0$ for every $r$, since otherwise the corresponding unit vector would be a nonzero proper integer kernel subvector of $(h_j,1,\ldots,1)$. Let $s_+$ and $s_-$ count the positive and negative values of $p^\top u^r$, respectively. Notice that $s_- \geq 1$ as $h_jk > 0$. Then Lemma~\ref{lem:product} gives 
\begin{equation}
\label{eq:two-row-count} (h_j+s_+)s_-\le\sum_{r:\,p^\top u^r<0}|p^\top u^r|. 
\end{equation}
We will use this inequality to prove parts~\ref{part:strict-column} and \ref{part:graver-coordinate} of the lemma.

For part~\ref{part:graver-coordinate}, we first show that $|p^\top u^r|\le\Delta(A_{-j})$ for every $r$. 
If $u^r$ uses a zero coefficient $q_i=0$, then $u^r=e_i$. Choose $\ell\ne j$ with $q_\ell\ne0$. Integrality gives 
\begin{equation}
    \label{eq:graver-q_i-0}
    |p^\top u^r|=|p_i|\le|p_iq_\ell|=\left|\det(a_i,a_\ell)\right|\le\Delta(A_{-j}).
\end{equation}
Otherwise, apply Theorem~\ref{thm:hilbert} to $u^r\in \mathcal{H}(q^\top)$ on the nonzero coefficients of $q$. For $q_i>0>q_\ell$, the corresponding normalized generator satisfies 
\begin{equation}
    \label{eq:two-row-normalized-generator}
    p^\top r_{i\ell}=p^\top((-q_\ell)e_i+q_i e_\ell)=-\det(a_i,a_\ell)
\end{equation}
by \eqref{eq:two-row-determinants}. The triangle inequality therefore gives $|p^\top u^r|\le\Delta(A_{-j})$. Using this estimate in \eqref{eq:two-row-count}, together with $s_- \geq 1$ and $s_+ \geq 0$, proves $h_j\le\Delta(A_{-j})$.

For part~\ref{part:strict-column}, suppose $h_e>0$ and $j\ne e$. At least one vector $u^r$ uses coordinate $e$. We first show that $|p^\top u^r| < \Delta$ for such an $r$. If $q_e=0$, we have $u^r = e_e$, and the same calculation as in \eqref{eq:graver-q_i-0} gives $|p^\top u^r|\le|\det(a_e,a_\ell)|<\Delta$.
If $q_e\ne0$, apply Theorem~\ref{thm:hilbert} to $u^r\in \mathcal{H}(q^\top)$. Since $u^r_e>0$, its convex representation has a positive coefficient on a normalized generator supported on a pair $\{i,e\}$. The same calculation as in \eqref{eq:two-row-normalized-generator} bounds the absolute first-row contribution of this generator by $|\det(a_i,a_e)|<\Delta$. Every other generator contributes at most $\Delta$ in absolute value. Since the coefficient sum is at most one, the triangle inequality yields $|p^\top u^r|<\Delta$.
The same estimates show that all $|p^\top u^r|$ are at most $\Delta$, and at least one is strictly smaller than $\Delta$. If $s_+>0$, then \eqref{eq:two-row-count} gives $h_j\le\Delta-s_+<\Delta$ as $s_-\geq 1$. If $s_+=0$, every $p^\top u^r$ is negative, so \eqref{eq:two-row-count} implies $h_j s_-\le\sum_{r:\,p^\top u^r<0}|p^\top u^r|<s_-\Delta$. 
This also gives $h_j<\Delta$ and proves part~\ref{part:strict-column}.
\end{proof}

\section{Proof of Two-Row Proximity}\label{sec:two-proof}
Recall from Section~\ref{sec:link-hilbert-bases} that $B$ is an optimal basis representing $x^*$. Set $\delta=\left|\det A_B\right|$. Then $1\le\delta\le\Delta$. Throughout this section, $z^*$ denotes the selected optimum of \eqref{eq:IP} satisfying \eqref{eq:M}. We first lift Theorem~\ref{thm:main} to a result for two rows in an extended Gomory group relaxation. The Gomory group relaxation~\cite{gomory1969polyhedra} releases the nonnegativity constraints on the basic variables $B$ in \eqref{eq:IP}. For this relaxation, tight proximity bounds for the $\ell_\infty$-norm are already available \cite[Theorem~1]{alievDistanceSparsityTransferenceVertices2021}. 
We investigate the extended Gomory group relaxation~\cite{wolsey1971extensions}, where, for $j\in B$, one releases only the nonnegativity constraint on coordinate $j$: 
\begin{equation}
\label{eq:IP-j} \mathrm{IP}_j(A,b,c)=\max\{c^\top z:Az=b,\ z\in\Z^n,\ z_{-j}\ge0\}.\tag{$\mathrm{IP}_j$} 
\end{equation}
The relaxation \eqref{eq:IP-j} is feasible because it contains every feasible point of \eqref{eq:IP}. Following the argument used for \eqref{eq:M}, we have $c^\top z=c^\top x^*+\bar c_N^\top z_N\le c^\top x^*$. Thus, since its integer hull is a rational polyhedron, \eqref{eq:IP-j} attains a finite optimum. 
For a feasible point $\widehat z$ of \eqref{eq:IP-j}, consider the minimality condition 
\begin{equation}
\tag{$\ast_j$}\label{eq:M-j} \nexists\,w\in\Z^n:Aw=b,\quad w_{-j}\ge0,\quad w_N\le\widehat z_N,\quad w_N\ne\widehat z_N. 
\end{equation}
An optimal solution $\widehat z^*$ of \eqref{eq:IP-j} minimizing $\sum_{i\in N}\widehat z_i^*$ satisfies \eqref{eq:M-j}, by the same reduced-cost argument used for \eqref{eq:M}. 
To prove a variant of Theorem~\ref{thm:main} for \eqref{eq:IP-j}, we fix $j\in B$ and use the standard form \eqref{eq:two-row-standard}. Each feasible solution to the equation $q^\top x_{-j} = \beta$ lifts to a solution of $Ax=b$, and each kernel vector $g_{-j}\in \ker(q^\top)$ lifts to a vector $g\in\ker(A)$. In both cases, the $j$-th coordinate is recovered by
\begin{equation}
\label{eq:two-row-projection} x_j=\frac{\gamma-p^\top x_{-j}}{k}, \qquad g_j=-\frac{p^\top g_{-j}}{k}.
\end{equation}
In particular, $x\in \Z^n$ if and only if $x_{-j}\in\Z^{n-1}$ and $\gamma - p^\top x_{-j} \equiv 0 \mod k$. Thus integer feasibility of \eqref{eq:IP-j} corresponds to nonnegative integer solutions of $q^\top x_{-j} = \beta$ satisfying $\gamma - p^\top x_{-j} \equiv 0 \mod k$. 
For three distinct indices, define the two-row normalized generator
\[
 r_{ijk}=\det(a_j,a_k)e_i-\det(a_i,a_k)e_j+\det(a_i,a_j)e_k.
\]
Let $\mathcal R_{B,j}$ consist of all nonzero vectors $\pm r_{ijh}$ with nonnegative $N$-part for distinct $i,h\ne j$. Notice that $\left|r_j\right|\le\Delta(A_{-j})$ for all $r\in\mathcal R_{B,j}$. For distinct $i,h\ne j$, one can calculate using \eqref{eq:two-row-determinants} that a one-row normalized generator $-q_h e_i+q_i e_h$ lifts to $(1/k) r_{jih}$ by \eqref{eq:two-row-projection}. Thus a nonnegative combination of such two-coordinate one-row generators with nonnegative $N$-parts lifts to a nonnegative combination of generators in $\mathcal R_{B,j}$ with each coefficient divided by $k$. Using these normalized generators, the next theorem gives a variant of Theorem~\ref{thm:main} that applies, in particular, to the optimal solution $\widehat z^*$ that satisfies \eqref{eq:M-j}.

\begin{theorem}\label{thm:extended-gomory}
Let $j\in B$. Every feasible point $\widehat z$ of \eqref{eq:IP-j} satisfying \eqref{eq:M-j} obeys 
\[
 \widehat z-x^*\in\left(1-\frac1\Delta\right)\conv(\{0\}\cup\mathcal R_{B,j}). 
\]
\end{theorem}

\begin{proof}
After a unimodular transformation we may use the standard form \eqref{eq:two-row-standard} with $B=\{j,\ell\}$ and $q_\ell=t>0$ such that $kt=\delta$ and $\beta=t x_\ell^*\ge0$. Consider the polyhedron defined by $q^\top x = \beta$ and $x\geq 0$. Mapping $\widehat z$ and $x^*$ into this polyhedron gives $u=\widehat z_{-j}$ and $u^* = (\beta/t)e_\ell$. The latter is a vertex of this polyhedron. 
Our goal is to bound the distance of $u^*$ to a nearby integer vector. Pick among all integer solutions of $q^\top w=\beta$, $w\ge0$, $w_N\le u_N$ one that minimizes $\sum_{i\in N}w_i$ and denote it by $v$. Notice that $q_i = 0$ implies $v_i = 0$ due to this minimality.
By Theorem~\ref{thm:main} and the optimality discussion surrounding \eqref{eq:M}, we can express $v-u^*\in \ker(q^\top)$ using normalized generators whose coefficient sum is at most $1-k/\Delta$, since $\Delta(q^\top)=\max_{i\ne j}|q_i|\le\Delta/k$ by \eqref{eq:two-row-determinants}.

Next we want to measure the distance from $u$ to $v$. To do so, decompose $u-v\in \ker(q^\top)$ into conformal Graver elements of $q^\top$: $u-v=g^1+\cdots+g^s$.
Since $u_N-v_N\geq 0$, every $g^r_N$ is nonnegative. Moreover, $g^r_N\ne0$: otherwise, $0=q^\top g^r=tg^r_\ell$ would imply $g^r=0$, since $t>0$.
We claim that no nonempty subset of $p^\top g^1,\ldots,p^\top g^s$ sums to zero modulo $k$. Suppose towards a contradiction that, without loss of generality, $p^\top(g^1+\ldots + g^r) \equiv 0 \mod k$. This implies that, after lifting $u - (g^1 + \ldots + g^r)$ back, its $j$-th coordinate is an integer by \eqref{eq:two-row-projection}. Using this, $g^i_N$ is nonzero, and $u - (g^1 + \ldots + g^r)$ stays componentwise between $u$ and $v$, implies that $u_N$ strictly decreases, contradicting \eqref{eq:M-j}. So no nonempty subset of $p^\top g^1,\ldots,p^\top g^s$ sums to zero modulo $k$. Lemma~\ref{lem:zero-sum} applied in $\Z/k\Z$ therefore gives $s\le k-1$.

Each $g^r$ is either a convex combination of zero and conformal one-row normalized generators by Theorem~\ref{thm:hilbert} or a signed unit vector at a zero coefficient. In the first case, as explained in the preceding discussion, extending a one-row normalized generator by its $j$-th coordinate according to \eqref{eq:two-row-projection} gives $1/k$ times a generator in $\mathcal R_{B,j}$. In the case when $g^r$ is a signed unit vector, we have $g^r\neq \pm e_{\ell}$ as $q_\ell=t>0$. Thus, as $g^r_N\geq 0$, there exists $i\neq \ell$ with $q_i = 0$ and $g^r = e_i = 1/t(-q_ie_{\ell} + t e_i)$. Its lift equals $-\frac1{kt}r_{ji\ell}$ by \eqref{eq:two-row-projection}.
Since $t\ge1$, this contributes weight at most $1/k$. So in both cases the lift of each $g^r$ has a representation using generators in $\mathcal R_{B,j}$ with total weight at most $1/k$. The lift of $v-u^*$ has a representation with total weight at most $(1-k/\Delta)/k$ by \eqref{eq:two-row-projection}. Lifting $u-u^*=(v-u^*)+\sum_{r=1}^s g^r$ back and estimating the total coefficient sum of the two-row normalized generators gives
\[
 \frac{1-k/\Delta+s}{k}\le\frac{1-k/\Delta+k-1}{k}=1-\frac1\Delta
\]
as $s\leq k-1$.
\end{proof}

The preceding theorem establishes a proximity result for the original integer program whenever $\widehat{z}_j\geq 0$. This is the case when $x^*_j$ is sufficiently large. 
\begin{corollary}\label{cor:transfer}
Let $z$ be an integer feasible point satisfying \eqref{eq:M}. If for some $j\in B$ 
\begin{equation}
\label{eq:transfer-condition} x_j^*\ge\left(1-\frac1\Delta\right)(\Delta(A_{-j})-1)\qquad\hbox{or}\qquad z_j\ge \Delta(A_{-j}), 
\end{equation}
 then $z$ satisfies \eqref{eq:M-j}.
\end{corollary}
\begin{proof}
Under the first condition choose $\widehat z$ feasible for \eqref{eq:IP-j} with $\widehat z_N\le z_N$ minimizing $\sum_{i\in N}\widehat z_i$. So $\widehat{z}$ satisfies \eqref{eq:M-j}.
Theorem~\ref{thm:extended-gomory} and the triangle inequality yield 
\begin{equation*}
    \widehat z_j\ge x_j^*-\left(1-\frac1\Delta\right) \Delta(A_{-j})\ge-\left(1-\frac1\Delta\right)>-1.
\end{equation*}
Its integrality makes $\widehat z_j$ nonnegative, so it is feasible for the original problem. Property~\eqref{eq:M} forces $\widehat z_N=z_N$, and linear independence then forces $\widehat z=z$.

Under the second condition, suppose $z$ does not satisfy \eqref{eq:M-j}, witnessed by $w$. Since $z$ satisfies \eqref{eq:M}, we must have $w_j<0$. Decompose $z-w\in \ker(A)$ into conformal Graver elements of $A$. Each has a nonzero and nonnegative nonbasic part, because a kernel vector only supported on $B$ is zero. Choose one of those Graver basis elements, say $g$. Outside $j$, subtracting $g$ from $z$ stays between $z$ and $w$, so remains nonnegative.
At $j$, $0\le g_j\le \Delta(A_{-j})\le z_j$ by Lemma~\ref{lem:two-row-bases}\ref{part:graver-coordinate}. Thus $z-g$ is feasible for the integer program and violates \eqref{eq:M}, a contradiction.
\end{proof}

We now combine Theorem~\ref{thm:extended-gomory}, Corollary~\ref{cor:transfer}, and Lemma~\ref{lem:augmented} to prove the two-row proximity result.
\begin{proof}[Proof of Theorem~\ref{thm:infinity}]
Fix an optimal LP vertex $x^*$ and the selected integer optimum $z^*$ satisfying \eqref{eq:M}. If $x^*$ is integral, this minimality forces $z^*=x^*$.
Otherwise, we bound the basic and nonbasic coordinates of their difference separately.

First fix $j\in B$. If either condition in \eqref{eq:transfer-condition} holds for $x^*$ and $z=z^*$, Corollary~\ref{cor:transfer} and Theorem~\ref{thm:extended-gomory} together with the triangle inequality give $|z_j^*-x_j^*|\le\left(1-\frac1\Delta\right)\Delta(A_{-j})\leq \Delta - 1$ as $\Delta(A_{-j})\le\Delta$. If neither condition holds, integrality gives $0\le z_j^*\le \Delta(A_{-j})-1\le \Delta - 1$, while $0\le x_j^*<\left(1-\frac1\Delta\right)(\Delta(A_{-j})-1)\le \Delta-1$. Both numbers then lie in $[0,\Delta-1]$ and thus $|z_j^*-x_j^*|\le \Delta-1$.

We use Lemma~\ref{lem:augmented} for the nonbasic variables. Assume $B=\{1,2\}$. Let $\bar x$ be obtained by either rounding $x^*_i$ up or down for $i=1,2$. In all four cases, Lemma~\ref{lem:augmented} implies $(|z^*-\bar x|,1)\in\mathcal H([AS,-(b-A\bar{x})])$ for some suitable sign matrix $S$. Consider first the rounded vector $\bar{x}$ obtained from rounding each $x^*_i$ to the nearest integer, i.e., $|\bar{x}_i - x^*_i|\leq 1/2$ for $i=1,2$. 
If $|\bar{x}_1 - x^*_1| + |\bar{x}_2 - x^*_2| < 1$, then 
\[
 \left|\det(b-A\bar{x},a_i)\right| \leq \left|\bar{x}_1 - x^*_1\right| \left|\det(a_1,a_i)\right| + \left|\bar{x}_2 - x^*_2\right|\left|\det(a_2,a_i)\right|<\Delta
\]
for all $i\in\lbrack n\rbrack$. So Lemma~\ref{lem:two-row-bases}\ref{part:strict-column} applied to $(|z^*-\bar x|,1)\in\mathcal H([AS,-(b-A\bar{x})])$ gives $z_i^*\leq \Delta - 1$ for all $i\in N$.

The remaining case is $|\bar{x}_i - x^*_i|= 1/2$ for $i=1,2$. Consider the rounded vector $\bar{x}$ that gives $b-A\bar{x} = (-a_1+a_2)/2$.
From the previous argument, we may assume that there exists a column $a_k$ such that $\left|\det(b-A\bar{x},a_k)\right|= \Delta$. So $a_k$ attains equality in the triangle inequality above and $\det(a_1,a_k)=-\det(a_2,a_k)=\pm \Delta$. 
Hence $a_k=\lambda(a_1+a_2)$ with $|\lambda|=\Delta/\delta$.
Since $a_k$ is a column of $A$, we get $|\det(a_1+a_2,a_i)|\le\delta$ for every column $a_i$. Now consider the rounded vector $\bar{y}$ that gives $b-A\bar{y}=(a_1+a_2)/2$. We just established that $\left|\det(b-A\bar{y},a_i)\right|\le \delta/2 < \Delta$ for all $i\in\lbrack n\rbrack$. Apply Lemma~\ref{lem:two-row-bases}\ref{part:strict-column} to $(|z^*-\bar y|,1)\in\mathcal H([AS',-(b-A\bar{y})])$, where $S'$ is the sign matrix associated with $z^*-\bar y$, to get $z_i^*\le\Delta-1$ for all $i\in N$.

Together with the basic-coordinate bounds, this gives $\|z^*-x^*\|_\infty\le\Delta-1$ for the same selected optimum. Taking the supremum over all admissible instances and optimal LP vertices gives $\pi^\infty(2,\Delta)\le\Delta-1$. Section~\ref{sec:lower-bound} supplies equality for every $\Delta\ge2$.
\end{proof}

\section{Lower Bound Examples}\label{sec:lower-bound}
For $\Delta = 1$, we may choose $z^*=x^*$ and all proximity bounds are zero and thus sharp. Fix $m$ and $\Delta\geq 2$. For a full-row-rank matrix $A$, both vectors $\pm r^I$ are called normalized generators, where
\[
 r^I=\sum_{\ell=0}^m(-1)^\ell\det A_{I\setminus\{i_\ell\}}\,e_{i_\ell},
 \qquad I=\{i_0<\cdots<i_m\},\quad \operatorname{rank}(A_I)=m.
\]
Consider $\mathcal{B}_{m,\Delta} = \lbrace -\Delta,\ldots,\Delta\rbrace^{m+1}$. Every normalized generator, restricted to its index set $I$, is a nonzero vector in $\mathcal{B}_{m,\Delta}$. We need the following converse statement.
\begin{lemma}\label{lem:cofactor-realization}
Let $v\in\mathcal B_{m,\Delta}$ satisfy $\gcd(v)=1$ and $\Vert v\Vert_\infty=\Delta$.
Then $v$ is a normalized generator of a full row rank matrix $A\in\Z^{m\times(m+1)}$ with $\Delta(A)=\Delta$ and $\gcd(A)=1$.
\end{lemma}
\begin{proof}
Since $\gcd(v) = 1$, it can be completed to a unimodular matrix $U$ whose last column is $v$.
Let $A$ consist of the first $m$ rows of $U^{-1}$. The cofactor formula for the inverse of $U^{-1}$ shows that both $v$ and $-v$ are normalized generators of $A$. Its maximal minors are, up to signs, the coordinates of $v$, so $\Delta(A)=\Vert v\Vert_\infty = \Delta$ and $\gcd(A)= \gcd(v)=1$.
\end{proof}

\begin{proposition}\label{prop:prescribed-generator-attainment}
Let $v\in\mathcal B_{m,\Delta}$ satisfy $\gcd(v)=1$ and $\Vert v\Vert_\infty=\Delta$.
There exist $A\in\Z^{m\times(m+1)}$ with $\gcd(A)=1$ and $\Delta(A)=\Delta$, integral $b,c$, and a normalized generator $r\in\{v,-v\}$ such that the unique optimal LP vertex $x^*$ and the unique optimal integer solution $z^*$ satisfy
\[
 z^*-x^*=\left(1-\frac1\Delta\right)r.
\]
\end{proposition}
\begin{proof}
Without loss of generality, let $|v_{m+1}|=\Delta$. Lemma~\ref{lem:cofactor-realization} gives a matrix $A\in \Z^{m\times (m+1)}$ with normalized generators $\pm v$, $\gcd(A)=1$, and $\Delta(A)=\Delta$. Pick the sign of $r\in\lbrace -v,v\rbrace$ such that $r_{m+1} = \Delta$. Define 
\begin{equation*}
    z_i^*=\max\{0,\lceil(1-1/\Delta)r_i\rceil\} \qquad x^* = z^* - (1-1/\Delta)r.
\end{equation*}
By construction, we have $x^*,z^*\geq 0$ and $x^*_{m+1} = 0$. Set $b = Az^* \in\Z^m$ and consider the feasible region $Ax=b, x\geq 0$. Notice that $x^*$ is a vertex of that region and all feasible solutions are of the form $x^* + \lambda r$ for $\lambda\geq 0$ as $\ker(A)$ is one-dimensional. Since $\gcd(r) = 1$ and $z^* = x^* + (1-1/\Delta)r$, there exists no integer solution of the form $x^* + t r$ for $0 \leq t < (1-1/\Delta)$. Picking $c = -e_{m+1}$ makes $x^*$ the unique optimal LP solution and $z^*$ the unique optimal integer solution.
\end{proof}

Together with the observation that no normalized generator chosen in Theorem~\ref{thm:main} can have greatest common divisor $\Delta$, Proposition~\ref{prop:prescribed-generator-attainment} shows that the scaling factor in Theorem~\ref{thm:main} is sharp. This translates to sharp proximity bounds for a fixed arbitrary norm $\Vert\cdot\Vert$ when $m = 1$ by picking a vector $v^*$ that attains $\max\{\Vert v\Vert: v\in\mathcal{B}_{1,\Delta}, ~\gcd(v)=1, ~\Vert v\Vert_\infty=\Delta\}$. We explicitly demonstrate this for $\ell_p$-norms to address Corollary~\ref{cor:allp}.

\begin{corollary}\label{cor:lp-lower-bound}
For all positive integers $m$ and $\Delta$,
\[
 \Delta-1\le\pi^\infty(m,\Delta),\quad \left(1-\frac1\Delta\right)\bigl(m\Delta^p+(\Delta-1)^p\bigr)^{1/p}\le\pi^p(m,\Delta)\quad(1\le p<\infty).
\]
In particular,
\[
 (m+1)\Delta-m-2+1/\Delta\le\pi^1(m,\Delta),\quad \left(1-\frac1\Delta\right)\sqrt{m\Delta^2+(\Delta-1)^2}\le\pi^2(m,\Delta).
\]
\end{corollary}
\begin{proof}
For $\Delta\ge2$, apply Proposition~\ref{prop:prescribed-generator-attainment} to 
\begin{equation*}
    v = (\Delta,\ldots,\Delta,\Delta-1)^\top.
\end{equation*}
The claim follows since $\Vert v\Vert_p^p=m\Delta^p+(\Delta-1)^p$ for $1\le p<\infty$ and $\Vert v\Vert_\infty=\Delta$.
\end{proof}
Setting $m=1$ gives the lower bounds for Corollary~\ref{cor:allp} and setting $m=2$ with $p=\infty$ gives the lower bound in Theorem~\ref{thm:infinity}.

\section*{Acknowledgments}
The authors used ChatGPT-Astra for proof exploration, manuscript review and literature cross-checking. ChatGPT-Astra suggested the new proof of Theorem~\ref{thm:hilbert}, the idea of lifting the Hilbert basis argument, and the examples in Appendix~\ref{app:conformal}. The normalized generator approach, the feasibility variants of Theorem~\ref{thm:main} and Corollary~\ref{cor:allp}, Theorem~\ref{thm:extended-gomory} and Corollary~\ref{cor:transfer}, the results in Section~\ref{sec:lower-bound}, and other intermediate results were obtained by the authors before using ChatGPT-Astra and were provided to it as input for the contributions described above. The authors wrote the manuscript and take full responsibility for its correctness, originality, and final content.

\bibliography{proximity}

\begin{thebibliography}{10}

\bibitem{aliev2026optimal}
Iskander Aliev, Gennadiy Averkov, W.~Jones, and Timm Oertel.
\newblock Optimal proximity bound and product function estimates in integer linear programming, 2026.
\newblock URL: \url{https://arxiv.org/abs/2606.13579}.

\bibitem{alievSparsityProximityTransference2025}
Iskander Aliev, Marcel Celaya, and Martin Henk.
\newblock Sparsity and proximity transference in integer programming.
\newblock {\em Mathematical Programming}, 216:29--48, 2026.
\newblock \href {https://doi.org/10.1007/s10107-024-02191-z} {\path{doi:10.1007/s10107-024-02191-z}}.

\bibitem{alievDistanceSparsityTransferenceVertices2021}
Iskander Aliev, Marcel Celaya, Martin Henk, and Aled Williams.
\newblock Distance-{Sparsity} {Transference} for {Vertices} of {Corner} {Polyhedra}.
\newblock {\em SIAM Journal on Optimization}, 31(1):200--216, 2021.
\newblock \href {https://doi.org/10.1137/20M1353228} {\path{doi:10.1137/20M1353228}}.

\bibitem{alievSupportIntegerOptimal2018}
Iskander Aliev, Jes\'{u}s~A. De~Loera, Friedrich Eisenbrand, Timm Oertel, and Robert Weismantel.
\newblock The {Support} of {Integer} {Optimal} {Solutions}.
\newblock {\em SIAM Journal on Optimization}, 28(3):2152--2157, 2018.
\newblock \href {https://doi.org/10.1137/17M1162792} {\path{doi:10.1137/17M1162792}}.

\bibitem{aliev2020distances}
Iskander Aliev, Martin Henk, and Timm Oertel.
\newblock Distances to lattice points in knapsack polyhedra.
\newblock {\em Mathematical Programming}, 182(1--2):175--198, 2020.
\newblock \href {https://doi.org/10.1007/s10107-019-01392-1} {\path{doi:10.1007/s10107-019-01392-1}}.

\bibitem{aprileIntegerProgramsNearly2024}
Manuel Aprile, Samuel Fiorini, Gwenaël Joret, Stefan Kober, Michał Seweryn, Stefan Weltge, and Yelena Yuditsky.
\newblock Integer {Programs} with {Nearly} {Totally} {Unimodular} {Matrices}: {The} {Cographic} {Case}.
\newblock {\em Mathematics of Operations Research}, 2026.
\newblock \href {https://doi.org/10.1287/moor.2024.0830} {\path{doi:10.1287/moor.2024.0830}}.

\bibitem{celaya2022proximity}
Marcel Celaya, Stefan Kuhlmann, Joseph Paat, and Robert Weismantel.
\newblock Proximity and flatness bounds for linear integer optimization.
\newblock {\em Mathematics of Operations Research}, 49(4):2446--2467, 2024.
\newblock \href {https://doi.org/10.1287/moor.2022.0335} {\path{doi:10.1287/moor.2022.0335}}.

\bibitem{cook1986sensitivity}
William Cook, A.~M.~H. Gerards, Alexander Schrijver, and {\'E}va Tardos.
\newblock Sensitivity theorems in integer linear programming.
\newblock {\em Mathematical Programming}, 34(3):251--264, 1986.
\newblock \href {https://doi.org/10.1007/BF01582230} {\path{doi:10.1007/BF01582230}}.

\bibitem{eisenbrand_proximity_2020}
Friedrich Eisenbrand and Robert Weismantel.
\newblock Proximity results and faster algorithms for integer programming using the {Steinitz} lemma.
\newblock {\em ACM Transactions on Algorithms}, 16(1):5:1--5:14, 2020.
\newblock \href {https://doi.org/10.1145/3340322} {\path{doi:10.1145/3340322}}.

\bibitem{gomory1969polyhedra}
Ralph~E. Gomory.
\newblock Some polyhedra related to combinatorial problems.
\newblock {\em Linear Algebra and its Applications}, 2(4):451--558, 1969.
\newblock \href {https://doi.org/10.1016/0024-3795(69)90017-2} {\path{doi:10.1016/0024-3795(69)90017-2}}.

\bibitem{gribanov$$Delta$$modularInteger2024}
Dmitry Gribanov, Ivan Shumilov, Dmitry Malyshev, and Panos Pardalos.
\newblock On {$\Delta$-modular} integer linear problems in the canonical form and equivalent problems.
\newblock {\em Journal of Global Optimization}, 88(3):591--651, 2024.
\newblock \href {https://doi.org/10.1007/s10898-022-01165-9} {\path{doi:10.1007/s10898-022-01165-9}}.

\bibitem{HenkKuhlmannWeismantel2022}
Martin Henk, Stefan Kuhlmann, and Robert Weismantel.
\newblock On {Lattice} {Width} of {Lattice}-{Free} {Polyhedra} and {Height} of {Hilbert} {Bases}.
\newblock {\em SIAM Journal on Discrete Mathematics}, 36(3):1918--1942, 2022.
\newblock \href {https://doi.org/10.1137/21M145416X} {\path{doi:10.1137/21M145416X}}.

\bibitem{HenkWeismantel2000}
Martin Henk and Robert Weismantel.
\newblock On minimal solutions of diophantine equations.
\newblock {\em Beiträge zur Algebra und Geometrie}, 41(1):49--55, 2000.

\bibitem{hochbaumConvexSeparableOptimization1990}
Dorit~S. Hochbaum and J.~George Shanthikumar.
\newblock Convex separable optimization is not much harder than linear optimization.
\newblock {\em Journal of the ACM}, 37(4):843--862, 1990.
\newblock \href {https://doi.org/10.1145/96559.96597} {\path{doi:10.1145/96559.96597}}.

\bibitem{koppe2026determining}
Matthias K{\"o}ppe, Moises Reyes~Rivas, and Luze Xu.
\newblock Determining sharp proximity bounds for low row rank and {$\Delta$}-modularity.
\newblock {\em INFORMS Journal on Optimization}, 2026.
\newblock \href {https://doi.org/10.1287/ijoo.2025.0068} {\path{doi:10.1287/ijoo.2025.0068}}.

\bibitem{kuhlmannThresholdPhenomenonShortest2026}
Stefan Kuhlmann and Robert Weismantel.
\newblock A {Threshold} {Number} for the {Shortest} {Vector} {Problem} in the {Infinity} {Norm}, 2026.
\newblock \href {https://doi.org/10.48550/arXiv.2508.02249} {\path{doi:10.48550/arXiv.2508.02249}}.

\bibitem{lee2020improving}
Jon Lee, Joseph Paat, Ingo Stallknecht, and Luze Xu.
\newblock Improving proximity bounds using sparsity.
\newblock In Mourad Ba{\"i}ou, Bernard Gendron, Oktay G{\"u}nl{\"u}k, and A.~Ridha Mahjoub, editors, {\em Combinatorial Optimization: ISCO 2020}, volume 12176 of {\em Lecture Notes in Computer Science}, pages 115--127, Cham, 2020. Springer.
\newblock \href {https://doi.org/10.1007/978-3-030-53262-8_10} {\path{doi:10.1007/978-3-030-53262-8_10}}.

\bibitem{schrijver1998theory}
Alexander Schrijver.
\newblock {\em Theory of linear and integer programming}.
\newblock John Wiley \& Sons, 1998.

\bibitem{sissokho2014note}
Papa~A. Sissokho.
\newblock A note on minimal zero-sum sequences over {$\mathbb Z$}.
\newblock {\em Acta Arithmetica}, 166(3):279--288, 2014.
\newblock URL: \url{http://eudml.org/doc/286310}.

\bibitem{sissokhoGeometryMinimalSolutions2021}
Papa~A. Sissokho.
\newblock Geometry of the {Minimal} {Solutions} of a {Linear} {Diophantine} {Equation}.
\newblock {\em SIAM Journal on Discrete Mathematics}, 35(3):1952--1963, 2021.
\newblock \href {https://doi.org/10.1137/21M1404417} {\path{doi:10.1137/21M1404417}}.

\bibitem{wolsey1971extensions}
Laurence~A Wolsey.
\newblock Extensions of the group theoretic approach in integer programming.
\newblock {\em Management Science}, 18(1):74--83, 1971.
\newblock \href {https://doi.org/10.1287/mnsc.18.1.74} {\path{doi:10.1287/mnsc.18.1.74}}.

\end{thebibliography}

\clearpage
\appendix

\section{Failure of Conformal Representations for $m=2$}\label{app:conformal}
Theorem~\ref{thm:main} and Theorem~\ref{thm:hilbert} admit representations using conformal normalized generators. This is not true anymore when $m=2$. We begin with the Hilbert basis statement.
\par\medskip\noindent\textit{Example.}\label{ex:hilbert-nonconformal}
Let
\[
 A=\begin{pmatrix}1&0&-1&1&-4\\0&1&-5&-1&3\end{pmatrix},
 \qquad h=(3,4,1,2,1)^\top.
\]
One can verify that indeed $h\in\mathcal H(A)$. Since every coordinate of $h$ is positive, the normalized generators conformal to $h$ are precisely the nonnegative ones:
\[
\begin{aligned}
 r_{123}&=(1,5,1,0,0)^\top,& r_{135}&=(23,0,3,0,5)^\top,\\
 r_{145}&=(1,0,0,3,1)^\top,& r_{234}&=(0,6,1,1,0)^\top,\\
 r_{245}&=(0,1,0,4,1)^\top.
\end{aligned}
\]
For $c = (1,3,0,5,0)^\top$, the linear functional $c^\top x$ is at most $23$ on these generators, whereas $c^\top h = 25$. This implies that $h$ is not in their convex hull.
In fact,
\[
 h=\frac3{23}r_{135}+\frac{14}{23}r_{234}+\frac8{23}r_{245},
\]
and the functional shows that the coefficient sum $25/23>1$ is the smallest possible among conformal nonnegative representations.
Allowing the nonconformal normalized generators
\[
 r_{235}=(0,23,4,0,-1)^\top,\qquad r_{345}=(0,0,-1,23,6)^\top
\]
gives
\[
 h=\frac3{23}r_{135}+\frac4{23}r_{235}+\frac2{23}r_{345},
\]
whose coefficient sum is $9/23<1$.
\par
Due to the close connection between proximity and Hilbert bases, it may not come as a surprise that an analogous statement exists for proximity.

\par\medskip\noindent\textit{Example.}
For every integer $q\ge3$, let $B=\{1,2\}$,
\[
 A=\begin{pmatrix}q&1&1&0\\1&q&0&1\end{pmatrix},\quad b=(q-1,q-1)^\top,\quad c=(0,0,-1,-1). 
\]
One has $\det A_B =\Delta=q^2-1$. The feasible region is bounded because the row sum is strictly positive on every column. Let $f=(q-1)/(q+1)$. The unique LP optimum is $x^*=(f,f,0,0)$.
If an integer basic coordinate were positive, the corresponding row would exceed its right-hand side. Therefore the unique integer feasible point is $z^*=(0,0,q-1,q-1)$. 
There are exactly four normalized generators with nonnegative nonbasic part: 
\[
 \begin{array}{ll} r_3=(-q,1,\Delta,0),&r_4=(1,-q,0,\Delta),\\ s_1=(-1,0,q,1),&s_2=(0,-1,1,q). \end{array} 
\]
 Only $s_1,s_2$ are conformal to $z^*-x^*=(-f,-f,q-1,q-1)$. Their coefficients are forced by the two basic entries to be $f,f$.
Thus the least total conformal weight is 
\[
 2f= 1- \frac{1}{q^2-1} + \frac{(q-2)^2}{q^2-1}>1-\frac{1}{\Delta}. 
\]
Allowing the nonconformal generators $r_3,r_4$ gives
\[
 z^*-x^*=\frac{r_3+r_4}{q+1},\qquad
 \frac2{q+1}\le1-\frac1\Delta.
\]
\par

\end{document}